\documentclass[12pt]{article}
\usepackage{amsmath}
\usepackage{amscd}
\usepackage{graphicx}
\usepackage{amssymb,amsfonts,amsthm,amscd,mathrsfs}
\usepackage{verbatim}
\usepackage[T2A]{fontenc}
\usepackage[utf8]{inputenc}
\usepackage[russian]{babel}
\usepackage{bm}
\usepackage{authblk}
\renewcommand{\b}{\beta}

\newcommand{\e}{\varepsilon}

\newcommand{\E}{\mathsf{E}}

\begin{document}


\title{  Local Power of Tests of Fit for  Normality of Autoregression
}
\author{M.V. Boldin
\footnote{Dept. Mech. Math., Moscow State Lomonosov  Univ., Moscow, Russia\\
 e-mail:boldin$_{-}$m@ hotmail.com}}

\date{ }
\maketitle



We consider a stationary $AR(p)$ model. The autoregression parameters are unknown as well as the distribution of innovations. Based on the residuals from the parameter estimates, an analog of empirical distribution function is defined and the tests of Kolmogorov's and $\omega^2$ type is constructed for testing hypotheses on the normality of innovations. We obtain the asymptotic power of these tests under local alternatives.


AMS {\bf 2010 Subject Classification:} Primary 62C10; secondary 62M10, 62G30.
\section{Введение }

 В этой работе мы изучаем асимптотическую локальную мощность критериев типа Колмогорова-Смирнова и Крамера-Мизеса-Смирнова омега-квадрат (далее кратко критериев Колмогорова и омега-квадрат) для проверки нормальности линейной стационарной авторегрессионной модели. Параметры модели неизвестны, так что инновации, порождающие авторегрессионную последовательность, ненаблюдаемы. Упомянутые тесты основываются на эмпирической функции распределения 
  остатков от подходящих оценок параметров. Эта функция строится по остаткам аналогично тому, как строится (гипотетически) эмпирическая функция распределения по самим инновациям.
 
 Подобные остаточные э.ф.р. изучались  во многих работах для различных авторегрессионных моделей. В частности, в \cite{Bold82} для стационарной $AR(p)$ модели с нулевым средним, в \cite{Bold83} для линейной регрессии с авторегрессионными ошибками, в \cite{Koul.Levent.1989} для взрывающейся авторегрессии, в \cite{Kreiss} для $AR(\infty)$, в \cite{Bold2000} для $ARCH(p)$ модели и т.д.
 
 В упомянутых работах на остаточной эмпирической функции распределения сновываются статистики типа  Колмогорова и омега-квадрат для проверки гипотез о распределении инноваций, найдены предельные распределения этих статистик при гипотезах. Мощность соответствующих тестов при локальных альтернативах до самого последнего времени была неизвестна и это было пробелом в упомянутых привлекательных фактах.
 
 (Строго говоря, термин ''локальная мощность'' подразумевает наличие регулярной последовательности стохастических экспериментов, см., например, \cite{LeCam}. Мы такую последовательность определять не будем. Под локальной мощностью мы понимаем далее мощность при близких альтернативах, сходящихся к гипотезе со скоростью $O(n^{1/2})$, $n$ -- объем данных.)
 
 В недавней работе \cite{BoldMMS2019} для линейной стационарной $AR(p)$ модели с нулевым средним удалось найти асимптотическую локальную мощность тестов типа Колмогорова и омега-квадрат для проверки простых гипотех о распределении инноваций. Оказалось, что эта асиптотическая мощность такая же, как в классической ситуации независимых одинаково распределенных наблюдений. Напомним, в этой классической ситуации локальные мощности критериев согласия были найдены впервые в \cite{Chib65}. 

 В настоящей заметке мы делаем следующий шаг в сравнении с \cite{BoldMMS2019}. Мы рассматриваем линейную стационарную $AR(p)$ авторегрессию с ненулевым средним и проверяем гипотезу о нормальности инноваций. Мы находим асимптотическую локальную мощность соответствующих тестов типа Колмогорова и омега-квадрат.
 
 Напомним, в линейной авторегрессии нормальность инноваций эквивалентна нормальности самой авторегрессионной последовательности и обеспечивает оптимальность общеупотребительных процедур наименьших квадратов. См., например, \cite{And}, гл. 5 и 6. Это объясняет теоретический и прикладной интерес к рассматриваемой далее задаче.

 Мы описываем модель и приводим некоторые предварительные результаты в Разделе 2,  формально ставим задачу и формулируем основной результат в Разделе 3. Доказательства весьма кропотливы и будут опубликованы в отдельной работе.

 \section{Описание модели и тестовых статистик}

 Рассмотрим стационарную AR(p) модель с ненулевым средним
$$
v_t = \b_1 v_{t-1} + \dots + \b_p v_{t-p} +\nu+ {\e}_t, \  t \in \mathbb{Z.} \eqno(1)
$$
В (1) $\{{\e}_t\}$ -- независимые одинаково распределенные случайные величины (н.о.р.сл.в.) с неизвестной функцией распределения (ф.р.) G(x);  $\E \e_1 = 0$,  $0<\E {\e}_1^2 < \infty$; $\bm{\b} = (\b_1, \dots, \b_p)^T$ - вектор неизвестных параметров, таких что корни соответствующего (1) характеристического уравнения по модулю меньше единицы; неизвестная константа $\nu\in \mathbb{R}$, порядок уравнения $p$ известен.\\
Пусть наблюдения $v_{1-p}, \dots, v_n$ - выборка из стационарного решения $\{{u}_t\}$ уравнения (1).

Построим тесты типа Колмогорова и омега-квадрат для проверки гипотезы
$$
H_0: \,G(x)\in\{\Phi(x/\sigma),\,\, \sigma>0\},\,\, \Phi(x) - \text{ф.р. стандартной гауссовской величины}.
$$
Сначала перепишем уравнение (1) в удобном для дальнейшего рассмотрения виде.\\
Введем константу $\mu$ соотношением $\nu=(1-\b_1-\ldots-\b_p)\mu$, тогда
$$
v_t-\mu = \b_1 (v_{t-1}-\mu) + \dots + \b_p (v_{t-p}-\mu) + {\e}_t. 
$$
Полагая $u_t:=v_t-\mu$, получаем:
$$
v_t=\mu+u_t,\eqno(2)
$$
$$
u_t = \b_1 u_{t-1} + \dots + \b_p u_{t-p} + {\e}_t, \  t \in \mathbb{Z.} \eqno(3)
$$
Соотношения (2) -- (3) эквивалентны (1), при этом (3) -- авторегрессия с нулевым средним. Вспоминая вид стационарного решения (3), см., например, \cite{And}, гл. 5, получаем строго стационарное решение (2) --(3) в виде
$$
v_t=\mu+\sum_{j\geq 0}\gamma_j\e_{t-j}.\eqno(4)
$$
Ряд в (4) сходится в $L^2$, последовательность $\{\gamma_j\}$ определяется рекуррентным соотношением
$$
\gamma_j = \b_1 \gamma_{t-1} + \dots + \b_p \gamma_{t-p},\,j\geq 1,\,\gamma_0=1,\,\gamma_j=0 \mbox{ для } \,j<0,
$$
$|\gamma_j|\leq c \lambda^j,\, 0<\lambda<1$, константы $c,\lambda$ от $j$ и вектора $\bm \b$ не зависят.

Соотношение (4) влечет: последовательность $\{v_t\}$ из (4) является гауссовской тогда и только тогда, когда верна $H_0$.

Вернемся к построению тестов для $H_0$. Возьмем оценкой $\mu$ в (2) эмпирическое среднее 
$
\overline \mu:=n^{-1}\sum_{t=1}^n v_t,
$ 
а оценками ненаблюдаемых $\{u_t\}$ величины
$$
\hat u_t:=v_t-\hat v,\quad t=1-p,\ldots,n.
$$
Пусть   $\bm{\hat \b_n} = (\hat \b_{1n}, \dots, \hat \b_{pn})^T$ - любая $n^{1/2}$ - состоятельная оценка вектора $\bm{\b}$, построенная по этим $\{\hat u_t\}$. Например, годится оценка наименьших квадратов (о.н.к.), поскольку она  асимптотически нормальна при сделанных предположениях, см., например, \cite{And}, гл. 5.
Величины 
$$
\hat \e_t = \hat u_t -  \hat \b_{1n} \hat u_{t-1} - \dots - \hat \b_{pn} \hat u_{t-p},\quad t = 1, \dots, n,
$$
называются отстатками, это оценки ненаблюдаемых инноваций $\{\e_t\}$, а функция
$$
\hat G_n(x) = n^{-1} \sum_{t=1}^n I(\hat \e_t \le x),\,x\in \mathbb R,
$$
называется остаточной эмпирической функцией распределения (о.э.ф.р). (Здесь и далее $I(\cdot)$ означает индикатор события.) $\hat G_n(x)$ есть подобие эмпирической функции распределения 
$$
G_n(x) = n^{-1} \sum_{t=1}^n I(\e_t \le x)
$$
ненаблюдаемых сл.в. $\e_{1}, \dots, \e_n$.

Асимптотические свойства $\hat G_n(x)$ характеризует следующее стохастическое разложение, см. \cite{Bold83}: 

если $G(x)$ дважды дифференцируема с $g(x) = G'(x)$ и $\sup_x|g'(x)| < \infty$, то 
$$
\sup_{x} |n^{1/2}[\hat G_n(x) - G_n(x+\overline \e)]|  \xrightarrow{P} 0, \quad n \to \infty, \quad \overline \e:=n^{-1}\sum_{t=1}^n \e_t\eqno(5)
$$
Неизвестную дисперсию инноваций при $H_0$ обозначим $\sigma_0^2$. Таким образом, гипотеза $H_0$ эквивалентна тому, что $G(x)=\Phi(x/\sigma_0)$ для некоторой неизвестной дисперсии $\sigma_0^2>0$.  
Оценкой $\sigma_0^2$ возьмем
$$
\hat s_n^2:=n^{-1}\sum _{t=1}^n \hat \e_t ^2.
$$
Введем остаточный эмпирический процесс (о.э.п.)
$$
\hat u_n(x):=n^{1/2}[\hat G_n(x) - \Phi(x/\hat s_n)].
$$
Разложение (5) позволяет указать слабый предел процесса $\hat u_n(\hat s_n \Phi^{-1}(t)),t\in [0,1]$, в 
пространстве Скорохода $D[0,1]$ функций без разрывов второго рода в предположении, что гипотеза $H_0$ верна.\\
А именно, пусть $u(t),t\in [0,1]$, будет гауссовский процесс со средним ноль и ковариацией
$$
\min\{t,s\}-ts-\varphi(\Phi^{-1}(t))\varphi(\Phi^{-1}(s))-\frac{1}{2}\Phi^{-1}(t)\varphi(\Phi^{-1}(t))\Phi^{-1}(s)\varphi(\Phi^{-1}(s)),
$$
$\varphi(x):=\Phi'(x),\,\Phi^{-1}(t)$ -- обратная к $\Phi(x)$ функция.
В \cite{Bold83} показано: при $H_0$
$$
 \hat u_n(\hat s_n \Phi^{-1}(t))\xrightarrow{D[0,1]} u(t),\quad n\to \infty.\eqno(6)
$$
Статистики типа Колмогорова и омега-квадрат для $H_0$ имеют вид:
$$
\hat D_n:=\sup_t|\hat u_n(\hat s_n \Phi^{-1}(t))|=\sup_x|\hat u_n(x)|,\quad 
$$
$$
\hat\omega_n^2=\int_0^1[\hat u_n(\hat s_n \Phi^{-1}(t))]^2dt=\int_{-\infty}^{\infty}[\hat u_n(x))]^2d \Phi(x/\hat s_n).
$$
В силу (6 ) при $H_0$
$$
P(\hat D_n \leq \lambda)\to P(\sup_t| u(t)|\leq \lambda):=\hat K(\lambda),
$$
$$
 P( \hat\omega_n^2 \leq \lambda)\to P(\int_0^1[ u(t)]^2 dt\leq \lambda):=\hat S(\lambda),\quad n\to \infty.
$$
Квантили  ф.р. $\hat K(\lambda)$,  оцененные методом Монте-Карло, приведены в таблицах \cite{PearsonHartly}, а квантили ф.р. $\hat S(\lambda)$  вычислены в \cite{Mart}. \\Таким образом,  статистики  $\hat D_n$ и $\hat\omega_n^2$ можно применять для проверки $H_0$ при больших $n$ так же, как обычные статистики Колмогорова и омега-квадрат. 

Упомянутые факты (в частности, соотношение (6)) были установлены только при гипотезе $H_0$, поведение о.э.п. $\hat u_n(\hat s_n \Phi^{-1}(t))$ и основанных на нем тестовых статистик при локальных альтернативах не рассматривалось. В следующем разделе мы устраним этот пробел.

\section{Постановка задачи и основной результат}

Цель настоящей работы -- установить слабый предел в $D[0,1]$ о.э.п. $\hat u_n((\hat s_n \Phi^{-1}(t))$ и основанных на нем статистик  
$\hat D_n$ и $\hat\omega_n^2$ при локальных альтернативах.\\
А именно, далее  предполагается, что $\{\e_t\}$ в (2) -- (3) -- н.о.р.сл.в. с ф.р. в виде смеси:
$$
G(x)=A_n(x):=(1-n^{-1/2})\Phi(x/\sigma_0)+n^{-1/2} H(x),\,H(x) - \mbox{ф.р.}\eqno (7)
$$

Предположение (7) будем понимать как локальную альтернативную к $H_0$ гипотезу и обозначать $H_{1n}$. Разумеется, при $H(x)=\Phi(x/\sigma_0)$альтернатива $H_{1n}$ и гипотеза $H_0$ совпадают. \\
При $H_{1n}$ вектор $\b$ остается постоянным и от $n$ не зависит. Авторегрессионное среднее $\mu$ от $n$ может зависеть и вот почему. Пусть среднее значение ф.р. $H(x)$ будет, скажем, $\nu_H \neq 0$. Тогда уравнения (2) -- (3) можно переписать в эквивалентном виде, заменив последовательность $\{\e_t\}$ на последовательность с нулевым средним $\{\e_t-n^{-1/2} \nu_H\}$, а $\mu$ на $\mu_n:=(1-\b_1-\ldots-\b_p)^{-1}(\nu+ n^{-1/2} \nu_H)$. Значит, без ограничения общности среднее инноваций при альтернативе $H_{1n}$ можно считать нулевым.

Сказанное делает естественным следующее предположение.

{\bf Условие (i).}
Случайная величина с функцией распределения  $H(x)$ имееет нулевое среднее и конечную дисперсию $\sigma_H^2$.

Отметим, что остатки $\{\hat \e_t\}$ , определенные в Разделе 2, ни при гипотезе ни при альтернативе от среднего авторегрессионной схемы $\mu$ вовсе не зависят. Потому свободны от этого среднего и процесс $\hat u_n(x)$ и статистики $\hat D_n,\,\hat \omega_n^2$.  

{\bf Условие (ii).}
Функция распределения $H(x)$ дифференцируема с производной, удовлетворяющей условию Липшица.

Условия (i) и (ii) позволяют использовать результаты \cite{BoldMMS2019} об асимптотическом поведении о.э.ф.р. при локальных альтернативах.

Пусть $\bm{\hat \b}_n$ будет любая  $n^{1/2}$ -- состоятельная при $H_{1n}$ 
 оценка $\bm{\b}$. Например, при $H_{1n}$ и единственном Условии (i) годится о.н.к. $\hat{\bm \b}_{n,LS}$, построенная по $\{\hat u_t\}$. Она остается асимптотически гауссовской при $H_{1n}$ с теми  же параметрами, как при $H_0$. А именно,
 $$
 n^{1/2}(\hat{\bm \b}_{n,LS}-\b) \xrightarrow{d} \mathbb N(\bm 0,\,\sigma_0^2 \bm K^{-1}),\quad n\to\infty.
 $$
 Здесь ковариационная матрица $\bm K=(k_{i,j})>0.\,i,j=1,\ldots,p\,\,$, $k_{i,j}=\E u_0^0 u_{i-j}^0$\,\,,$\{u_{t}^0\}$ -- стационарное решение уравнения (3) с $\{\e_t\}$, имеющими ф.р. $\Phi(x/\sigma_0)$.

 Доказательство этого факта такое же, как в авторегрессии с нулевым средним, см. \cite{BoldMMS2019}, а потому опущено.\\

 Вот наш основной результат.
\newtheorem{Theorem}{Теорема} 
\begin{Theorem}
Пусть верна альтернатива $H_{1n}$. Пусть выполнены Условия (i) и (ii). Положим 
$$
\delta(t):=H(\sigma_0 \Phi^{-1}(t))-t +\frac{1}{2}\Phi^{-1}(t)\varphi(\Phi^{-1}(t))(\frac{\sigma_H^2}{\sigma_0^2}-1),\quad t\in [0,1].
$$
Тогда 
$$
\hat {u}_n(\hat s_n \Phi^{-1}(t))=n^{1/2}[\hat G_n(\hat s_n \Phi^{-1}(t)) -  t] \xrightarrow{D[0,1]} u(t)+\delta(t),\quad n\to \infty.
$$
\end{Theorem}
Теорема 1 прямо влечет
\newtheorem{Corollary}{Следствие} 
\begin{Corollary}
В условиях теоремы 1 справедлива сходимость по распределению:  
$$
\hat D_n\xrightarrow{d} \sup_t|u(t)+\delta(t)|,\,  \quad    \hat \omega^2_n \xrightarrow{d} \int_0^1[ u(t)+\delta(t)]^2dt,                \quad n\to \infty.
$$
\end{Corollary}

Выделим в нашем рассмотрении важный частный случай н.о.р. наблюдений. Он получается при $\bm \b=0$, и в этом случае последовательность $\{v_t\}$  имеет вид
$$
v_t=\nu+\e_t, \quad t\in \mathbb Z.
$$
Полагая $\hat {\bm \b_n}=0$, получим остатки $\hat \e_t=\e_t-\overline \e,\,t=1,\ldots, n$. Для таких остатков 
$$
\hat G_n(x)=G_n(x+\overline \e), \hat s_n^2=s_n^2:=n^{-1}\sum_{t=1}^n(\e_t-\overline \e)^2,
$$
$$
\hat {u}_n(\hat s_n \Phi^{-1}(t))=u_n(t):=n^{1/2}[ G_n(\overline \e+ s_n \Phi^{-1}(t)) -  t],
$$
$$
\hat D_n=D_n:=\sup_t|u_n(t)|,\quad \hat \omega_n^2= \omega_n^2:=\int_{0}^1[u_n(t)]^2dt.
$$
Статистики $D_n$ и $\omega_n^2$ -- обычные статистики Колмогорова и омега-квадрат с оцененными параметрами, предназначеные для проверки нормальности н.о.р.сл.в. $v_1,\ldots,v_n$.

Теорема 1 и Следствие 1 прямо влекут 
\begin{Corollary}
Пусть верна альтернатива $H_{1n}$  и выполнены Условия (i) и (ii). Тогда
$$
 {u}_n(t)\xrightarrow{D[0,1]} u(t)+\delta(t),
$$
  $$
 D_n\xrightarrow{d} \sup_t|u(t)+\delta(t)|,\,  \quad     \omega^2_n \xrightarrow{d} \int_0^1[ u(t)+\delta(t)]^2dt,                \quad n\to \infty.
$$
\end{Corollary}
Заметим, что Условие (ii) в Следствии 2 излишне. Анализируя доказательство Теоремы 1 в Разделе 4 легко увидеть, что в частном случае $\bm \b=\hat{\bm \b}=0$ достаточно потребовать 
непрерывности $H(x)$.

В силу Следствия 2 слабые пределы при $H_{1n}$ процесса $\hat {u}_n(\hat s_n \Phi^{-1}(t))$ и cтатистик $\hat D_n, \hat \omega^2_n$ в общем случае те же, что у $ {u}_n(t), D_n$ и $ \omega^2_n$, построенных по н.о.р.сл.в. 
Весьма примечательный факт.

Отметим, что утверждения Следствия 2 известны давно, их можно получить из общих результатов \cite{Durbin1973}, где для н.о.р. наблюдений изучалась слабая сходимость в $D[0,1]$ эмпирического процесса с оцененными параметрами.

\end{document}